\documentclass[journal,onecolumn]{IEEEtran}

\usepackage{amssymb}%,bbm}
\usepackage{amsmath,amsthm}
\usepackage{subfigure}
\usepackage{cite}
\usepackage{graphicx}
\usepackage{graphics}
\usepackage{color}
\usepackage{xspace}
\usepackage{bbm}
\usepackage{psfrag}
\usepackage{algorithmicx}
\usepackage{algorithm}
\usepackage{algpseudocode}
\usepackage{multirow}
\usepackage{array}
\usepackage{url}
\usepackage[normalem]{ulem}
\usepackage{setspace}
\usepackage{mathtools}

\DeclarePairedDelimiter\floor{\lfloor}{\rfloor}

\usepackage{graphicx}
\usepackage{floatflt,setspace}
\usepackage{algcompatible}
\algnewcommand\INPUT{\item[\textbf{Input:}]}%
\algnewcommand\OUTPUT{\item[\textbf{Output:}]}%

\usepackage[mathscr]{eucal}
\usepackage{amsbsy}
\usepackage{bm}
\usepackage{fixltx2e}
\MakeRobust{\overrightarrow}

\include{macrofile}

\newcommand{\argmin}{{\text{argmin}}}

\title{Tensor Completion by Alternating Minimization under the Tensor Train (TT) Model}
\author{Wenqi Wang, Vaneet Aggarwal, and Shuchin Aeron\thanks{W. Wang and V. Aggarwal are with the School of Industrial Engineering, Purdue University, West Lafayette, IN, 47907, email: \{wang2041,vaneet\}@purdue.edu. S. Aeron is with the Dept. of Electrical and Computer Engineering, Tufts University, Medford, MA 02155, email: shuchin@ece.tufts.edu }}

\newtheorem{definition}{Definition}
\newtheorem{lemma}{Lemma}

\begin{document}
\pagenumbering{gobble}

\maketitle

% *******************************************************  Introduction  ********************************************************
\begin{abstract}
Using the matrix product state (MPS) representation of tensor train decompositions, in this paper we propose a tensor completion algorithm which alternates over the matrices (tensors) in the MPS representation. This development is motivated in part by the success of matrix completion algorithms which alternate over the (low-rank) factors. We comment on the computational complexity of the proposed algorithm and numerically compare it with existing methods employing low rank tensor train approximation for data completion as well as several other recently proposed methods. We show that our method is superior to existing ones for a variety of real settings.  
%Tensor train decomposition has shown to be effective in investigating the embedded hierarchy structure and compressing high dimensional tensor data.
%In this paper, a tensor completion algorithm under Tensor Train (TT) model is formulated.
%We further propose a Tensor Completion algorithm by Alternative Minimization under TT model (TCAM-TT) to solve the optimization problem, which solves the minimization problem within each alternative step by a computation effective least square method.  
%Numerical results show that TCAM-TT algorithm is effective in image and video completion under high missing ratios. 
\end{abstract}

\section{Introduction}
Tensor decompositions for representing and storing data have recently become very popular due to their effectiveness in effectively compressing data for statistical signal processing, see \cite{SIREV,CichockiMPCZZL14,tensorface1} for some of the applications. In this paper we focus on Tensor Train (TT) decomposition \cite{oseledets2011tensor} and in particular its relation to Matrix Product States (MPS) \cite{Orus:2013kga} representation for completing data from missing entries. In this context our algorithm is motivated by recent work in matrix completion where under a suitable initialization an alternating minimization algorithm \cite{jain2013low, Hardt13a} over the low rank factors is able to accurately predict the missing data. 

Tensor completion based on TT decompositions have been recently considered in \cite{phien2016efficient}. These approaches do not explicitly exploit the MPS representation of the TT format and therefore are not able to take the full advantage of this structured decomposition. Further our algorithm works by choosing a spectral initialization using just the available data, which results in reducing the number of iterations required for convergence for the proposed method. The proposed algorithm gives the detailed steps for solving the least square with respect to one of the tensor in the MPS representation. 

The rest of the paper is organized as follows. In section \ref{sec:2} we introduce the basic notation and preliminaries on the TT decomposition. In section \ref{sec:3} we outline the problem statement and propose the main algorithm in section \ref{HLSA_section}. Section \ref{compl} describes the computational complexity of the proposed algorithm. Following that we test the algorithm extensively against competing methods on a number of real and synthetic data experiments in section \ref{sec:5}. Finally we provide conclusion and future research directions in section \ref{sec:6}.

%Our novelty: 
%
%(1) Proposed a tensor completion algorithm under tensor train models, which is also potentially adapted in tensor ring models \cite{zhao2016tensor}. 
%
%(2) Proposed an alternative minimization method that updates each optimal solution in alternative minimization step by computational effective least square method .
%
%(3) Different from SiLRTC-TT algorithm as proposed in \cite{phien2016efficient}, which solves tensor train completion by weighted summarization of each unfolding completion, alternative minimization algorithm solves the tensor completion by taking orders to updated each tensor train factorization term.

% ****************************************************  Problem Set-up  ********************************************************
\section{Notation \& Preliminaries}
\label{sec:2}
In this paper, vector and matrices are represented by bold face lower case letters $({\bf x,y,z,\cdots})$ and bold face capital letters $({\bf X, Y, Z,\cdots})$ respectively. 
A tensor with order more than two is represented by calligraphic letters $(\bf \mathscr{X}, \mathscr{Y}, \mathscr{Z})$. 
For example, a $n^\text{th}$ order tensor is represented by ${\bf \mathscr{X}} \in \mathbb{R}^{I_1 \times I_2 \times \cdots \times I_n}$, where $I_{i: i=1,2,\cdots, n}$ is the tensor dimension along mode $i$. 
The tensor dimension along mode $i$ may be an expression, where the expression inside $()$ is evaluated as a scalar, e.g. $\mathscr{X}\in \mathbb{R}^{(I_1I_2) \times (I_3I_4)\times (I_5I_6)}$ represents a 3-mode tensor where dimensions along each mode is $I_1I_2$, $I_3I_4$, and $I_5I_6$ respectively. 

An entry inside a tensor $\mathscr{X}$ is represented as $\mathscr{X}(i_1, i_2,\cdots, i_n)$, where $i_{k: k=1,2,.., n}$ is the location index along the $k^{\text{th}}$ mode. 
A colon is applied to represent all the elements of a mode in a tensor,  e.g. $\mathscr{X}(:, i_2,\cdots, i_n)$ represents the fiber along mode $1$ and $\mathscr{X}[:, :, i_3, i_4,\cdots, i_n]$ represents the slice along mode $1$ and mode $2$ and so forth. 

Product notation $\otimes$ represents matrix Kronecker product and $\circ$ represents Hadamard product. Similar to Hadamard product under matrices case, Hadamard product between tensors is the entry-wise product of the two tensors.  
$\text{vec}(\cdot)$ represents the vectorization of the tensor in the argument. The vectorization is carried out lexicographically over the index set, stacking the elements on top of each other in that order. 
Frobenius norm of a tensor is the same as the vector $\ell_2$ norm of the corresponding tensor after vectorization, e.g. $\|\mathscr{X}\|_F =\|\text{vec}(\mathscr{X})\|_{\ell_2}$.

\begin{figure}[h!]
\includegraphics[height = 2.75in, width = 6in]{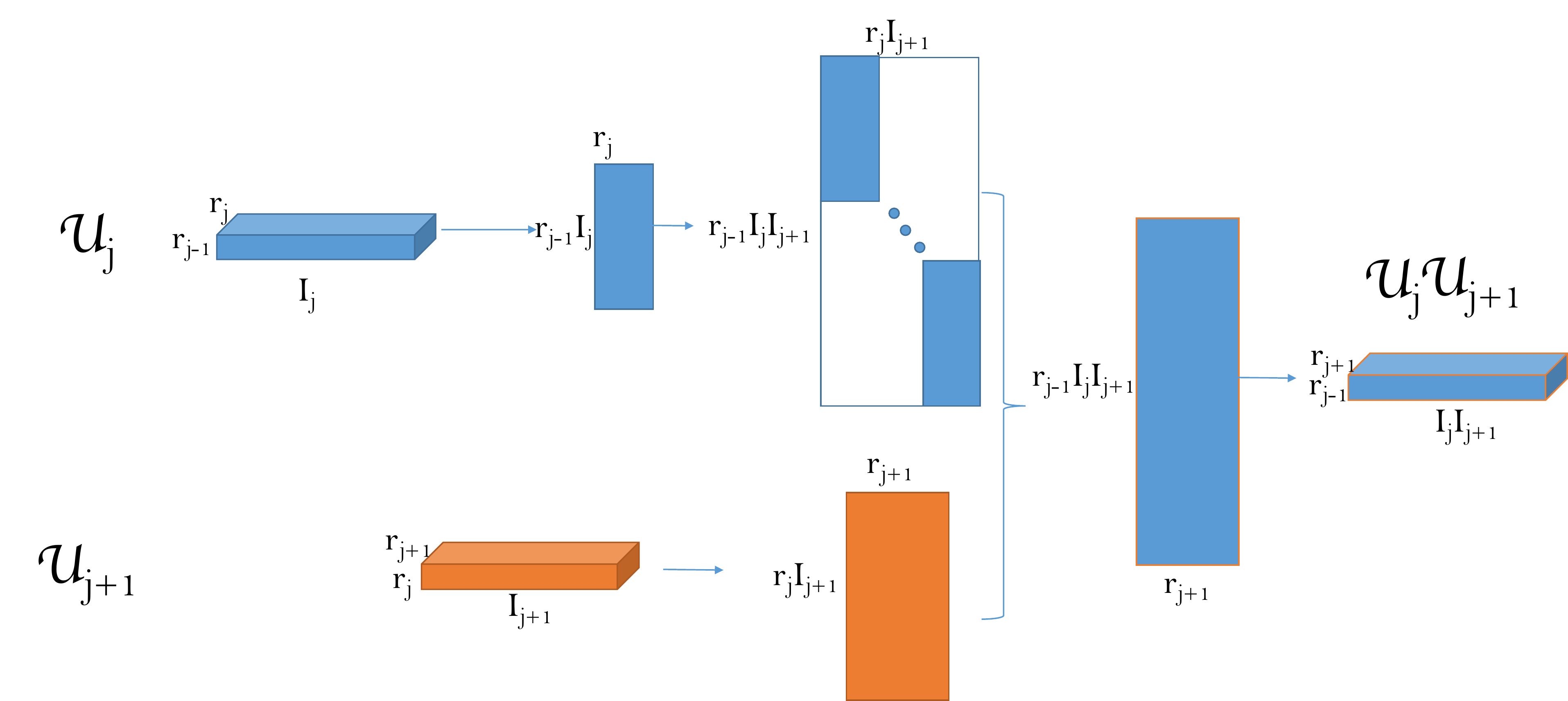}
\centering
\vspace{0.1in}
\caption{\normalsize{Tensor Connect Product for the $\mathscr{U}_j$ and $\mathscr{U}_{j+1}$}. 
For tensor connect product with more than $2$ tensors, connect product the first two tensors and take the connect product tensor to connect product with the third tensor. 
%\textcolor{red}{Wenqi: You need to tell how this extends to connect product of more than two tensors. Are they done in some order? Or is it order independent like multiplying a set of matrices together?}
}
\label{ConProd}
\end{figure}

We first introduce three commonly  used tensor unfolding operations namely, Tensor Mode-k Unfolding, Tensor Mode Matricization (TMM), Left-unfolding, and Right-unfolding as they will be intensively used in this paper.
% ================================   Definition 1   =======================================================================
%\textcolor{red}{What is the difference between the Mode-k unfolding and Tensor mode matricization? I though they are the same? Are we using some different operations here for the two. We need to be a little careful in using the language here.}

{\definition(Tensor Mode-k Unfolding)\\
The mode-$k$ unfolding matrix of a $n^{\text{th}}$ order tensor $\mathscr{X} \in \mathbb{R}^{I_1 \times \cdots \times I_n}$, denoted as ${\bf X}_{[k]} \in 
\mathbb{R}^{I_k \times (I_{k+1} \cdots I_nI_1 \cdots I_{k-1} )}$, such that
\begin{equation}
\begin{split}
&\mathscr{X}(i_1, i_2,\cdots, i_n)\\
=&{\bf X}_{[k]}(i_k, i_{k+1} + (i_{k+2}-1)I_{k+1} +\cdots \\
+&(i_{k-1}-1)\prod_{t=k+1,\cdots, n,1, \cdots k-1} I_t).
\end{split}
\end{equation}
}

% ================================   Definition 2   =======================================================================
{\definition (Tensor Mode Matricization (TMM)) \label{def2}\\
Let  $\mathscr{X}\in\mathbb{R}^{I_1 \times \cdots \times I_n}$ be a $n^{\text{th}}$ order tensor, the tensor mode matrization along the $k^{\text{th}}$ mode, denoted as $\mathscr{X}_{\floor{k}} \in \mathbb{R}^{(\prod_{t=1}^k I_t) \times (\prod_{t=k+1}^n I_t)}$, is a matrix where
\begin{equation}
\begin{split}
&\mathscr{X}(i_1,\cdots, i_n) \\
= &\mathscr{X}_{\floor{i}}( i_1 + (i_2-1)I_1 +\cdots+ (i_k-1) \prod_{t=1}^{k-1}I_t, \\
&i_{k+1}+(i_{k+2}-1)I_{k+1} + \cdots + (i_n-1) \prod_{t=k+1}^{n-1}I_t).
\end{split}
\end{equation}
\enddefinition}

% ================================  Definition 3 ===============================================================================
{\definition (Left-unfolding \& Right-unfolding \cite{holtz2012manifolds})\\
Let $\mathscr{U} \in \mathbb{R}^{r_0 \times I_1 \times r_1}$ be a 3rd-{order} tensor, the left-unfolding of $\mathscr{U}$, denoted as ${\bf L}(\mathscr{U})$, satisfies the following property 
\begin{equation}
{\bf L}(\mathscr{U}) \in\mathbb{R}^{(r_0I_1) \times r_1}, 
\end{equation}
and
\begin{equation}
\begin{split}
 {\bf L}(\mathscr{U})(k_0+(i_1-1)r_0, k_1) = \mathscr{U}(k_0, i_1, k_1)\\
 \forall k_0 \in [1, r_0], k_1\in[1, r_1], i_1 \in [1, I_1].
 \end{split}
\end{equation}
And let ${\bf L}^{-1}$ be the reverse operation of ${\bf L}$, which reshapes a $\mathbb{R}^{(r_0I_1)\times r_1}$ matrix to a $\mathbb{R}^{r_0\times I_1 \times r_1}$ tensor.

Similarly, the right unfolding of $\mathscr{U}$, denoted as ${\bf R}(\mathscr{U})$, satisfies the following property
\begin{equation}
{\bf R}(\mathscr{U}) \in\mathbb{R}^{r_0 \times (I_1r_1)}, 
\end{equation}
and
\begin{equation}
\begin{split}
 {\bf R}(\mathscr{U})(k_0, i_1 + (k_1-1)I_1) = \mathscr{U}(k_0, i_1, k_1)\\
 \forall k_0 \in [1, r_0], k_1\in[1, r_1], i_1 \in [1, I_1],
 \end{split}
\end{equation}
and ${\bf R}^{-1}$ is the reverse operation of ${\bf R}$, which reshapes a $\mathbb{R}^{r_0\times (I_1r_1)}$ matrix to a $\mathbb{R}^{r_0\times I_1 \times r_1}$ tensor.
\enddefinition}

%========================================================================================================
Tensor train decomposition \cite{oseledets2011tensor, holtz2012manifolds} is a tensor factorization method that any elements inside a tensor $\mathscr{X} \in \mathbb{R}^{I_1 \times\cdots\times I_n}$, denoted as $\mathscr{X}(i_1, i_2,\cdots, i_n)$, is represented by
\begin{align}
&\mathscr{X}(i_1,\cdots, i_n) \nonumber \\ &= {\bf U}_1(i_1,:) \mathscr{U}_2(:, i_2, :) \cdots \mathscr{U}_{n-1}(:, i_{n-1}, :){\bf U}_n(:, i_n),
\end{align}
where ${\bf U}_1\in\mathbb{R}^{I_1 \times r_1}$, ${\bf U}_n\in\mathbb{R}^{r_{n-1} \times I_n}$ are the boundary matrices and $\mathscr{U}_i \in \mathbb{R}^{r_{i-1} \times I_i \times r_i}, i=2,\cdots, n-1$ are middle decomposed tensors. 

A more general format of tensor train decomposition regards ${\bf U}_1 \in \mathbb{R}^{I_1 \times r_1}$ as a tensor $\mathscr{U}_1 \in \mathbb{R}^{r_0 \times I_1 \times r_1}, r_0=1$ and ${\bf U}_n \in \mathbb{R}^{r_{n-1} \times I_n}$ as a tensor $\mathscr{U}_n \in \mathbb{R}^{r_{n-1} \times I_n \times r_n}, r_n=1$, which gives the general tenor train decomposition format
\begin{equation}\label{eq: decomp1}
\mathscr{X}(i_1,\cdots, i_n) = \mathscr{U}_1(:, i_1,:) \cdots \mathscr{U}_{n-1}(:, i_{n-1}, :)\mathscr{U}_n(:, i_n, :),
\end{equation}
where $\mathscr{U}_i \in \mathbb{R}^{r_{i-1} \times I_i \times r_i}, i=1,\cdots, n$ and $r_0 = r_n=1$. The set of scalars,  $[r_0, r_1, \cdots, r_{n-1}, r_n]$, is defined as the tensor train rank (TT-Rank).

Since $\mathscr{X}(i_1,\cdots, i_n)$ is a scalar, \eqref{eq: decomp1} is equivalent to
 \begin{equation}\label{eq: tt-decom}
 \begin{split}
\mathscr{X}(i_1,\cdots, i_n) = \text{Trace}(&\mathscr{U}_1(:, i_1,:) \mathscr{U}_2(:, i_2, :) \cdots \\
&\mathscr{U}_{n-1}(:, i_{n-1}, :)\mathscr{U}_n(:, i_n, :)),
\end{split}
\end{equation}
which enables the cyclic permutations property ({See Definition 5 below}) that is  used  intensively in this paper. Before defining this property, we define Tensor Connect Product, that describes the product of a sequence of $3$rd-order tensor.

% ================================   Definition 4   =======================================================================
{\definition(Tensor Connect Product)
Let $\mathscr{U}_i \in \mathbb{R}^{r_{i-1} \times I_i \times r_i}, i=1,\cdots, n$ be $n$ $3$rd-order tensor, the tensor connect product is defined as,
\begin{equation}
\mathscr{U} =\mathscr{U}_1 \cdots \mathscr{U}_n \in \mathbb{R}^{r_0 \times (I_1\cdots I_n) \times r_n}.
\end{equation}
{and is shown in Fig \ref{ConProd}, where the tensor $\mathscr{U}_j$ is left-unfolded, denoted as ${\bf L}(\mathscr{U}_j)$ and the tensor $\mathscr{U}_{j+1}$ is left-unfolded, denoted as ${\bf L}(\mathscr{U}_{j+1})$. Then 
\begin{equation}
\begin{split}
\mathscr{U}_j\mathscr{U}_{j+1} & \in \mathbb{R}^{r_{j-1} \times (I_jI_{j+1}) \times r_{j+1}}\\
&= {\bf L}^{-1}({\bf I}^{(I_{j+1})} \otimes {\bf L}(\mathscr{U}_j) \times {\bf L}(\mathscr{U}_{j+1}))
\end{split}
\end{equation}}
}

%\begin{figure}[h]
%\includegraphics [trim=0.5in 0.6in 0.5in 0.5in, keepaspectratio, width=0.45\textwidth] {ConProd}
%\centering
%\vspace{0.1in}
%\caption{ Tensor Connect Product for the $\mathscr{U}_j$ and $\mathscr{U}_{j+1}$}
%\label{ConProd}
%\end{figure}

%Tensor connect product rule is shown in Fig \ref{ConProd}, where the tensor $\mathscr{U}_j$ is left-unfolded, denoted as ${\bf L}(\mathscr{U}_j)$ and the tensor $\mathscr{U}_{j+1}$ is right-unfolded, denoted as ${\bf R}(\mathscr{U}_{j+1})$. Then 
%\begin{equation}
%\begin{split}
%\mathscr{U}_j\mathscr{U}_{j+1} & \in \mathbb{R}^{r_{j-1} \times (I_jI_{j+1}) \times r_{j+1}}\\
%&= {\bf L}^{-1}({\bf I}^{(I_{j+1})} \otimes {\bf L}(\mathscr{U}_j) \times {\bf R}(\mathscr{U}_{j+1}))
%\end{split}
%\end{equation}

Let $f$ be a function applied on $\mathscr{U}$ such that $\mathscr{X} =f(\mathscr{U}) \in \mathbb{R}^{I_1\times \cdots \times I_n}$ satisfies \eqref{eq: tt-decom}, then $f$ is the function that reshapes vector $(I_1\cdots I_n) \times 1$ to tensor $I_1\times I_2 \times \cdots \times I_n$ after applied trace operation on each slice the $\mathscr{U}$ along mode-$2$, denoted as 
\begin{equation}
\mathscr{X} = f(\mathscr{U}),
\end{equation}
or equivalently
\begin{equation} \label{eq: eq6}
\mathscr{X} = f(\mathscr{U}_1 \cdots \mathscr{U}_n).
\end{equation}

Tensor connect product gives the product rule for the production between $3^{\text{rd}}$-order tensor, just like the matrix product as for $2^{\text{nd}}$ order tensor. We further note that tensor connect product is the same as matrix product for $2^{\text{nd}}$-order tensor.
{\lemma (Matrix Product) Tensor connect product is applied to matrix product. Let $M_1 \in \mathbb{R}^{I_1 \times r_1}$ and ${\bf M}_2 \in \mathbb{R}^{r_1 \times I_2}$ be any two matrix. Without loss of generality, we regard that ${\bf M}_1$ as a tensor $\mathscr{M}_1\in\mathbb{R}^{1 \times I_1 \times r_1}$ and ${\bf M}_2$ as $\mathscr{M}_2 \in \mathbb{R}^{r_1 \times I_2 \times 1}$ , then tensor connect product gives the vectorized solution of matrix production
\begin{equation}
\mathscr{M}_1 \mathscr{M}_2 \in \mathbb{R}^{1 \times (I_1I_2) \times 1}= \text{vec}({\bf M}_1{\bf M}_2).
\end{equation}
}
\proof
Proof is in Appendix \ref{proof0}.
\endproof

% ================================   Definition 5    =======================================================================
Similar to matrix transpose, which can be regarded as an operation that cyclic swaps the two modes for a $2^{\text{nd}}$ order tensor, 
%Similar to matrix transpose for $2^{\text{nd}}$ order tensor that cyclic swap the two mode of a tensor \textcolor{red}{$\leftarrow$ You mean cyclic swap of the indices?}, 
we defined Tensor Permutation to describe the cyclic-wise swap of tensor mode for high order tensor.
{\definition (Tensor Permutation) For any order-$d$ tensor $\mathscr{X} \in \mathbb{R}^{I_1 \times \cdots \times I_d}$, the $i^{\text{th}} $ tensor train permutation is defined as $\mathscr{X}^{P_i} \in \mathbb{R}^{I_i \times I_{i+1} \times \cdots \times I_n \times I_1 \times I_2 \times \cdots \times I_{i-1}}$  such that
\begin{equation}\label{eq: TensorPermutation}
\mathscr{X}^{P_i}(j_i,\cdots, j_n,j_1,\cdots, j_{i-1}) = \mathscr{X}(j_1,\cdots,j_n ), \forall_i, j_i\in [1,I_i].
\end{equation}
\enddefinition}

{We note the following result}.

\begin{lemma} \label{lemma1} $\mathscr{X}^{P_i} =f( \mathscr{U}_i \mathscr{U}_{i+1} \cdots \mathscr{U}_n \mathscr{U}_1 \cdots \mathscr{U}_{i-1})$.
\end{lemma} 
\proof
Proof is in Appendix \ref{proof1}
\endproof

With this background and basic constructs we now outline the main problem set-up.

% ==================================   Problem Setup   ==================================================================
\section{Problem Setup}
\label{sec:3}
Given a tensor $\mathscr{X} \in \mathbb{R}^{I_1 \times\cdots \times I_n}$ that is partially observed at locations $\Omega$, let $\mathscr{P}_\Omega \in \mathbb{R}^{I_1 \times\cdots \times I_n}$ be the corresponding binary tensor in which $1$ represents an observed entry and $0$ represents a missing entry.
The problem is to find a low tensor train rank (TT-Rank) approximation of the tensor $\mathscr{X}$, denoted as $\mathscr{W}$,  such that the recovered tensor $\mathscr{W}$ matches $\mathscr{X}$ at $\mathscr{P}_\Omega$. This problem is referred as the tensor completion problem under tensor train model,  which is equivalent to the following problem
\begin{align}\label{eq: T0}
\min_{\mathscr{W}:\mathscr{W}  \text{ satisfies TT-Rank } {\bf r}} \| \mathscr{P}_\Omega \circ ( \mathscr{W} -\mathscr{X}) \|_F^2 .
\end{align}

Using the factored form of TT representation, i.e. using equation \eqref{eq: decomp1}, the above optimization problem is equivalent to solving the following problem,
\begin{equation}\label{eq: T1}
\min _{\mathscr{U}_{i:i=1,2,\cdots,n}} \| \mathscr{P}_\Omega \circ (f(\mathscr{U}_1 \mathscr{U}_2\cdots \mathscr{U}_n)) -\mathscr{X}_\Omega) \|_F^2, 
\end{equation}
where the constraint that $\mathscr{W}$ is a low TT rank tensor is captured via $\mathscr{W} =\mathscr{U_1} \cdots \mathscr{U}_n$.

To solve this problem, We propose an algorithm referred to as Tensor Completion Algorithm by Alternating Minimization under the Tensor Train model, for short TCAM-TT, that solves the completion problem {in two steps, 
\begin{itemize}
\item Choosing an initial starting point by using Tensor Train Approximation (TTA) using the missing data only. This initialization algorithm is detailed in section \ref{TTA_section}.
\item Updating the solution by applying Hierarchical Alternating Least Square (HALS) that alternatively (in a cyclic order) estimates a factor say $\mathscr{U}_{i}$ keeping the other factors fixed. This algorithm is detailed in Section \ref{HLSA_section}.
\end{itemize}}

%*************************************    Tensor Train Decomposition    **********************************************************
\section{Tensor Train Approximation  (TTA) } \label{TTA_section}
%\textcolor{red}{Vaneet: Are you taking care of this section? This needs some work.}

For a given tensor $\mathscr{X}$, we wish to find the tensor $\mathscr{W}$ of TT-rank ${\bf r}$ that best approximates $\mathscr{X}$. Thus, we want to solve the problem given by

\begin{align}\label{approx}
\min_{\mathscr{W}:\mathscr{W}  \text{ satisfies TT-Rank } {\bf r}} \|  ( \mathscr{W} -\mathscr{X}) \|_F^2 .
\end{align}

Rather than solving the problem \eqref{approx} exactly, we give a heuristic algorithm to solve this problem. This is used as an initialization for the tensor completion problem, where the best approximation of zero-filled tensor will be used as an initialization. To avoid the computation complexity of \eqref{approx}, Algorithm 1 is used for the approximation. This algorithm gives the decomposition terms  $\mathscr{U}_i$ of the approximate solution $\mathscr{W}$.

%used a heuristic, this is because the approximation will be used as an initialization to our algorithm and thus it is enough to use a descent heuristic rather than optimally solving the exact approximation. 

%where SVD is performed rather than ?

%Tensor train decomposition factorizes a tensor following a hierarchy order where singular value decomposition(SVD) is applied to threshold the size of ranks. 
%Our proposed tensor train decomposition approximation is a modified version of tensor train decomposition as proposed in \cite{kammerer2011stability}. We reproduce the algorithm to fit our proposed tensor train format, as shown in Algorithm \ref{Algo_Decomposition}.

\begin{algorithm}
   \caption{Tensor Train Approximation}
   \label{Algo_Decomposition}
   \begin{algorithmic}[1]
   \INPUT  Tensor $\mathscr{X} \in \mathbb{R}^{I_1\times I_2\times \cdots \times I_n}$, TT-rank $r_{i: i=1,2,\cdots, n-1}$, $r_0 = r_n = 1$
    \OUTPUT Tensor train decomposition $\mathscr{U}_i \in \mathbb{R}^{r_{i-1}\times I_i \times r_{i}}, {i: i=1,2,\cdots, n}$
    \STATE \text{\bf Tensor Unfolding: } Apply tensor mode matrization for $\mathscr{X}$ along mode $1$ to get matrix ${\bf X}_1 =\mathscr{X}_{\floor{1}} \in \mathbb{R}^{I_1 \times (I_2 I_3 \cdots I_n)}$
    \STATE Apply SVD and threshold the number of singular values to be $r_1$ such that ${\bf X}_1 = {\bf U}_1 {\bf S}_1 {\bf V}_1^\top, 
    {\bf U}_1 \in \mathbb{R}^{I_1 \times r_1}, {\bf S}_1 \in \mathbb{R}^{r_1 \times r_1}, {\bf V}_1 \in \mathbb{R}^{r_1 \times (I_2 I_3 \cdots I_n)}$.  
    Note that ${\bf U}_1 = {\bf L}(\mathscr{U}_1)$, thus reshape ${\bf U}_1$ to $\mathbb{R}^{r_0 \times I_1 \times r_1}$ to recover $\mathscr{U}_1$ and let ${\bf M}_1 ={\bf S}_1 {\bf V}_1^\top \in \mathbb{R}^{r_1 \times ({I_2 I_3 \cdots I_n}) }$
	\FOR{ $i=2 $ to $n-1$}
	\STATE Reshape ${\bf M}_{i-1} \in \mathbb{R}^{r_{i-1} \times (I_i I_{i+1} \cdots  I_n)}$ to ${\bf X}_i \in \mathbb{R}^{(r_{i-1} I_i) \times (I_{i+1} I_{i+2} \cdots I_{n})}$ 
	\STATE Compute SVD and threshold the number of singular values to be $r_i$ such that ${\bf X}_{i} = {\bf U}_i {\bf S}_i {\bf V}_i^\top, 
	{\bf U}_i \in \mathbb{R}^{(r_{i-1}  I_i) \times r_i}, {\bf S}_i \in \mathbb{R}^{r_i \times r_i}, {\bf V} \in \mathbb{R}^{r_i \times (I_{i+1} I_{i+2} \cdots I_n)} $. 
	Note that ${\bf U}_i ={\bf L}(\mathscr{U}_i)$, thus reshape ${\bf U}_i \in \mathbb{R}^{(r_{i-1}I_i) \times r_i}$ to $\mathbb{R}^{r_{i-1} \times I_i \times r_i}$ to get $\mathscr{U}_i$ and set ${\bf M}_i ={\bf S}_i {\bf V}_i^\top \in \mathbb{R}^{r_i \times (I_{i+1} I_{i+2} \cdots  I_n)}$
	\ENDFOR
	\STATE Reshape ${\bf M}_{n-1} \in \mathbb{R}^{r_{n-1} \times I_n}$ to $\mathbb{R}^{r_{n-1} \times I_n \times r_n}$ to get $\mathscr{U}_n$
	\STATE Return $\mathscr{U}_1,\cdots, \mathscr{U}_n$
\end{algorithmic}
\label{main_algo}
\end{algorithm}

The proposed algorithm is a modified version of the tensor train decomposition as proposed in \cite{kammerer2011stability}. In the tensor train decomposition algorithm of \cite{kammerer2011stability}, the tensor is exactly  TT-Rank  ${\bf r}$. However, in our problem, the tensor $\mathscr{X}$ is not necessarily a TT-Rank  ${\bf r}$ tensor.
Thus, the singular value decomposition (SVD) is performed in different modes and thresholded to obtain the approximate TT-Rank  ${\bf r}$ tensor.

%%******************************************* Tensor Completion ****************************************************************
\section{Hierarchical Alternating Least Square (HALS)} \label{HLSA_section}
The proposed Tensor Completion method by Alternating Minimization under Tensor Train model (TCAM-TT) solves \eqref{eq: T1} by taking orders to solve the following problem %(\textcolor{red}{Wenqi: Can we use some other variable, other than $\mathscr{Y}$ below? We are using this for describing the problem.})
\begin{equation} \label{eq: T2}
\mathscr{U}_i = \argmin_{\mathscr{Y}} \| \mathscr{P}_\Omega \circ f(\mathscr{U}_1 \cdots \mathscr{U}_{i-1}\mathscr{Y}\mathscr{U}_{i+1}\cdots \mathscr{U}_n) -\mathscr{X}_\Omega) \|_F^2
\end{equation}

We further note that $\mathscr{U}_i$ in \eqref{eq: T2} can be solved
%\textcolor{red}{Wenqi: You mean Least Squares? Else this does not make sense to me. Can you explain why it is ALS for the LS problem?$\rightarrow$} 
%Alternative Least Square (ALS) method and 
only considering the following optimization problem
\begin{equation}\label{eq: T3}
\mathscr{U}_1 = \argmin_{\mathscr{Y} \in \mathbb{R}^{r_0 \times I_1 \times r_1}} \| \mathscr{P}_\Omega \circ f(\mathscr{Y} \mathscr{U}_2\cdots \mathscr{U}_n) -\mathscr{X}_\Omega) \|_F^2.
\end{equation}

%%%%%%%%%%%%%%%%%%%%%%%%%%%%%%%  Lemma 2  %%%%%%%%%%%%%%%%%%%%%%%%%%%%%%%%%%%%%%%%%%
{\lemma \label{lemma2}
When $i \neq 1$, solving 
\begin{equation} \label{eq: tran1}
\mathscr{U}_i = \argmin_{\mathscr{Y}} \| \mathscr{P}_\Omega \circ f(\mathscr{U}_1 \cdots \mathscr{U}_{i-1}\mathscr{Y}\mathscr{U}_{i+1}\cdots \mathscr{U}_n) -\mathscr{X}_\Omega) \|_F^2
\end{equation}
is equivalent to
\begin{equation}\label{eq: transform}
\mathscr{U}_i = \argmin_{\mathscr{Y}} \| \mathscr{P}^{P_i}_\Omega \circ f(\mathscr{Y} \mathscr{U}_{i+1}\cdots \mathscr{U}_n\mathscr{U}_{1}\cdots\mathscr{U}_{i-1}) -\mathscr{X}^{P_i}_\Omega  \|_F^2.
\end{equation} 
Since the format of \eqref{eq: transform} is exactly the same as \eqref{eq: T3}, thus solving $\mathscr{U}_i$ is equivalent to solving $\mathscr{U}_1$.
}
\proof
Proof is in Appendix \ref{proof2}
\endproof

%Based on Lemma \ref{lemma2}, without loss of generality, we only consider the scenario of updating $\mathscr{U}_1$ in equation \eqref{eq: T3}. The updating step for $\mathscr{U}_i$  follows the same logic besides pre-processing tensors by applying tensor permutation on $\mathscr{P}_\Omega$ and $\mathscr{X}_\Omega$ and changing orders of $U_i$, as described in \eqref{eq: transform}. 

\begin{figure}[h!]
\includegraphics [height = 3.25 in, width = 6.5in]{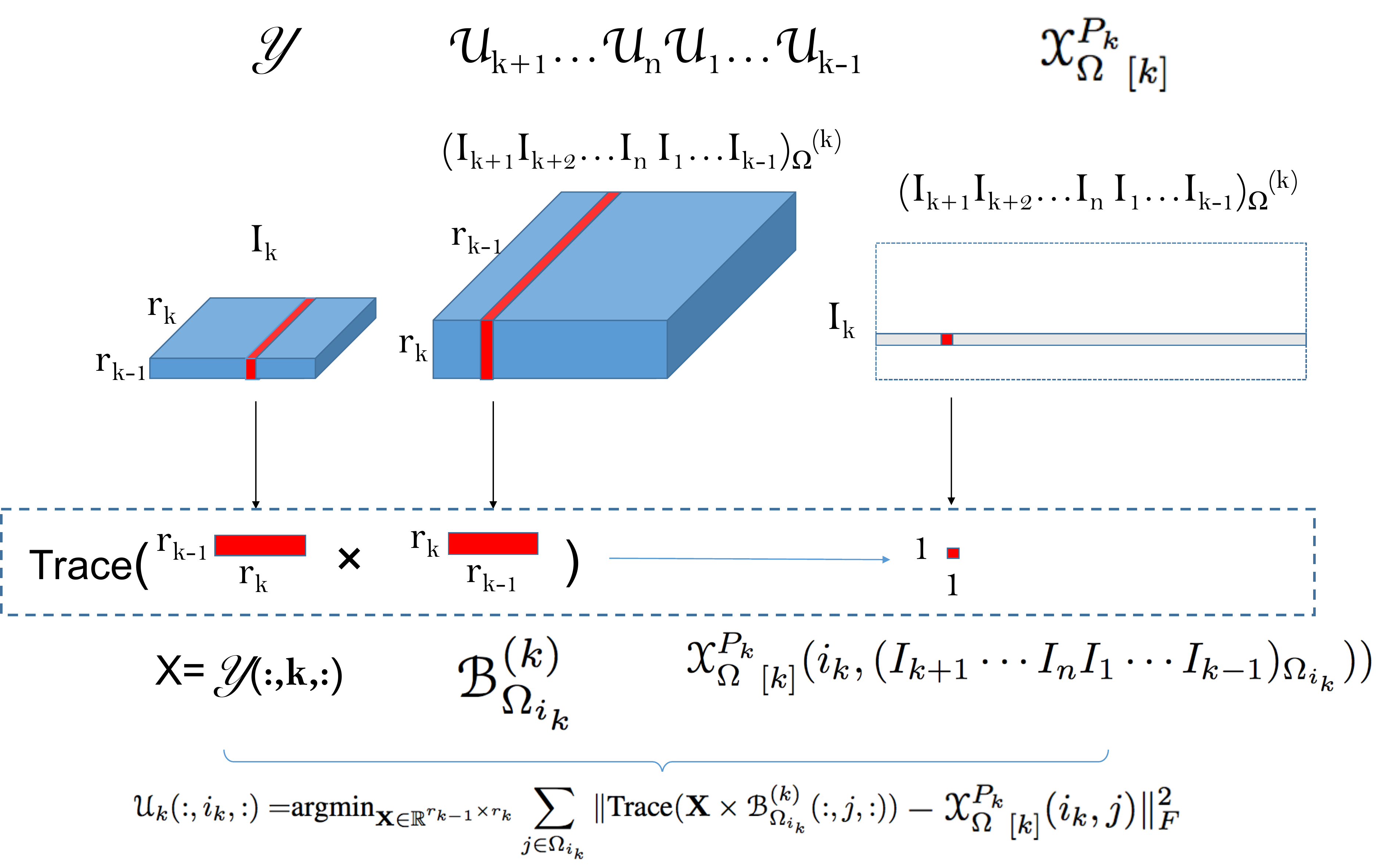}
\centering
\vspace{0.1in}
\caption{\normalsize{Updating $\mathscr{U}_k$ by Hierarchical Alternating Least Square (HALS)}}
\label{LS_U1}
\end{figure}

Now we consider solving $\mathscr{U}_k$ without loss of generality. Based on Lemma \ref{lemma2}, we need to solve the following problem
\begin{equation}\label{eq: problemk}
\mathscr{U}_k = \argmin_{\mathscr{Y}} \| \mathscr{P}^{P_k}_\Omega \circ f(\mathscr{Y} \mathscr{U}_{k+1}\cdots \mathscr{U}_n\mathscr{U}_{1}\cdots\mathscr{U}_{k-1}) -\mathscr{X}^{P_k}_\Omega  \|_F^2.
\end{equation} 

We further apply tensor mode-$k$ unfolding, which gives the equivalent problem
\begin{equation}\label{eq: problem1}
\mathscr{U}_k = \argmin_{\mathscr{Y}} \| {\mathscr{P}^{P_k}_\Omega}_{[k]} \circ {f(\mathscr{Y} \mathscr{U}_{k+1}\cdots \mathscr{U}_n\mathscr{U}_{1}\cdots\mathscr{U}_{k-1})}_{[k]} -{\mathscr{X}^{P_k}_\Omega}_{[k]}  \|_F^2.
\end{equation} 
where ${\mathscr{P}^{P_k}_{\Omega}}_{[k]}$, $ {f(\mathscr{Y} \mathscr{U}_{k+1}\cdots \mathscr{U}_n\mathscr{U}_{1}\cdots\mathscr{U}_{k-1})}_{[k]}$ and ${\mathscr{X}^{P_k}_\Omega}_{[k]} $ are matrices with dimension $\mathbb{R}^{I_k \times (I_{k+1}\cdots I_nI_1 \cdots I_{k-1})}$.

The trick in solving \eqref{eq: problem1} is that each slice of tensor $\mathscr{Y}$, denoted as $\mathscr{Y}(:, i_k, :), i_k=1,\cdots, I_k$ which corresponds to each row of ${\mathscr{P}^{P_k}_{\Omega}}_{[k]}$, $ {f(\mathscr{Y} \mathscr{U}_{k+1}\cdots \mathscr{U}_n\mathscr{U}_{1}\cdots\mathscr{U}_{k-1})}_{[k]}$ and ${\mathscr{X}^{P_k}_\Omega}_{[k]} $,   can be solved independently, thus equation \eqref{eq: T3} can be solved by solving $I_k$ equivalent subproblems 
\begin{equation}\label{eq: T4}
\begin{split}
&\mathscr{U}_k(:, i_k, :) \\
= &\argmin_{{\bf X} \in \mathbb{R}^{r_{k-1} \times 1 \times r_k}} \|{\mathscr{P}^{P_k}_\Omega}_{[k]}(i_k, :) \circ f({\bf X}\mathscr{U}_{k+1}\cdots \mathscr{U}_{k-1}) \\
&- {\mathscr{X}^{P_k}_\Omega}_{[k]}(i_k, :)\|_F^2
\end{split}
\end{equation}

%\begin{equation}\label{eq: T4}
%\begin{split}
%\mathscr{U}_k(:, i_k, :) = \argmin_{{\bf X}} \| &\mathscr{P}^{P_k}_\Omega \circ f({\bf X}\mathscr{U}_{k+1}\cdots \mathscr{U}_{k-1}) \\
%-&\mathscr{X}^{P_k}_\Omega(i_k,:,\cdots,:) \|_F^2.
%\end{split}
%\end{equation}

%======================================================   Slicing Part   =============================================================
As shown in Fig \ref{LS_U1}.

Let $\mathscr{B}^{(k)} = \mathscr{U}_{k+1} \cdots \mathscr{U}_n \mathscr{U}_1 \cdots \mathscr{U}_{k-1} \in \mathbb{R}^{r_k \times (I_{k+1} \cdots I_nI_1 \cdots I_{k-1}) \times r_{k-1}}$. 
Let $\mathscr{B}^{(k)}_{\Omega_{i_k}} \in \mathbb{R}^{r_k \times (I_{k+1} \cdots I_nI_1 \cdots I_{k-1})_{\Omega_{i_k}} \times r_{k-1}}$ be the components in $\mathscr{B}^{(k)}$ 
such that ${\mathscr{P}^{P_k}_\Omega}_{[k]}(i_k, (I_{k+1} \cdots I_nI_1 \cdots I_{k-1})_{\Omega_{i_k}})$ are observed.

Thus equation \eqref{eq: T4}   is equivalent to
\begin{equation}\label{eq: T5}
\begin{split}
\mathscr{U}_k(:, i_k,:) =&\argmin_{\mathscr{Z}} \| f(\mathscr{Z}\mathscr{B}^{(k)}_{\Omega_{i_k}}) \\
&- {\mathscr{X}^{P_k}_\Omega}_{[k]}(i_k, (I_{k+1} \cdots I_nI_1 \cdots I_{k-1})_{\Omega_{i_k}})) \|_F^2
\end{split}
\end{equation}
where 
$\mathscr{Z} \in \mathbb{R}^{r_{k-1} \times 1 \times r_k}$, 
$\mathscr{B}^{(k)}_{\Omega_{i_k}} \in \mathbb{R}^{r_k \times (I_{k+1} \cdots I_{k-1})_{\Omega_{i_k}} \times r_{k-1}}$, 
${\mathscr{X}^{P_k}_\Omega}_{[k]}(i_k, (I_{k+1} \cdots I_nI_1 \cdots I_{k-1})_{\Omega_{i_k}})) \in \mathbb{R}^{1 \times (I_{k+1} \cdots I_{k-1})_{\Omega_{i_k}} }$

 We regard $\mathscr{Z} \in \mathbb{R}^{r_{k-1} \times 1 \times r_k}$ as a matrix ${\bf X} \in \mathbb{R}^{r_{k-1} \times r_k}$. 
 Since the Frobenius norm of a vector in \eqref{eq: T5} is equivalent to entry-wise square summation of all entries, we rewrite \eqref{eq: T5} as 
\begin{equation}\label{eq: T5-1}
\begin{split}
\mathscr{U}_k(:, i_k,:) =&\argmin_{{\bf X}\in \mathbb{R}^{r_{k-1} \times r_k}} \sum_{j \in \Omega_{i_k}} \| \text{Trace}({\bf X} \times \mathscr{B}^{(k)}_{\Omega_{i_k}}(:, j, :))\\
 &-  {\mathscr{X}^{P_k}_\Omega}_{[k]}(i_k, j)\|_F^2
\end{split}
\end{equation}
where $\times$ is the matrix product.

{\lemma \label{lemma3} Let ${\bf A} \in \mathbb{R}^{r_1 \times r_2}$ and ${\bf B} \in \mathbb{R}^{r_2 \times r_1}$ be any two matrices, then 
\begin{equation}
\begin{split}
\text{Trace}({\bf A} \times {\bf B}) &= vec({\bf A})^\top vec({\bf B}^\top) =  vec({\bf B}^\top)^\top vec({\bf A})
\end{split}
\end{equation}
}
\proof 
Proof is in Appdendix \ref{proof3}.
\endproof

Based on Lemma \ref{lemma3}, \eqref{eq: T5-1} becomes
\begin{equation}\label{eq: T5-2}
\begin{split}
\mathscr{U}_k(:, i_k,:) =\argmin_{\bf X} \sum_{j \in \Omega^{(k)}_{i_k}} \|&vec((\mathscr{B}^{(k)}_{\Omega_{i_k}}(:, j, :))^\top)^\top vec({\bf X}) \\
&-  {\mathscr{X}^{P_k}_\Omega}_{[k]}(i_k, j)\|_F^2
\end{split}
\end{equation}

Then the problem for solving $\mathscr{U}_k[:, i_k, :]$ becomes a least square problem. Solving $I_k$ least square problem would give the optimal solution for $\mathscr{U}_k$.
Since each $\mathscr{U}_{i:i=1,\cdots,n}$ can solved by a least square method, tensor completion under tensor train model can be solved by taking orders to update $U_{i: i=1,\cdots, n}$ until convergence.

\begin{algorithm}
   \caption{TCAM-TT Algorithm}
   \begin{algorithmic}[1]
   \INPUT  Zero-filled Tensor $\mathscr{X}_\Omega \in \mathbb{R}^{I_1\times I_2\times ... \times I_n}$, binary observation index tensor $\mathscr{P}_\Omega \in \mathbb{R}^{I_1\times I_2\times ... \times I_n}$, TT-rank $r_{i: i=1,2,..., n-1}$, thresholding parameter $tot$, maximum iteration $maxiter$ 
    \OUTPUT Recovered tensor $\mathscr{X}_R$
    \STATE \text{\bf Tensor Train Decomposition Approximation}  Apply tensor train decomposition approximation in Algorithm 1 on $\mathscr{X}_\Omega$ to initialize the decomposition terms $\mathscr{U}_{i:i=1,..., n}^{(0)}$. Set iteration parameter $\ell = 0$. 
	\WHILE{ $\ell \leq maxiter$ }
	\STATE $\ell = \ell +1$
	\FOR {$i=1$ to $n$}
	\STATE \text{\bf Solve by Least Square Method} 
	${\mathscr{U}_i}^{(\ell)} = \argmin_\mathscr{U} \|\mathscr{P}_\Omega \circ (\mathscr{U} \mathscr{U}_{i+1}^{(\ell-1)}...\mathscr{U}_{n}^{(\ell-1)}\mathscr{U}_{1}^{(\ell)}...\mathscr{U}_{i-1}^{(\ell)}    - \mathscr{X})\|_F^2$ 
	\ENDFOR	
	\IF {$\sum_i \frac{ \|\mathscr{U}_i^{(\ell+1)} -\mathscr{U}_i^{(\ell)} \|_F}{ \| {\mathscr{U}}^{(\ell)}_i\|_F} \leq tot$}
	\STATE Break
	\ENDIF
	\ENDWHILE
	\STATE {\bf Recover Completed Tensors} $\mathscr{X}_R = \mathscr{U}_1^{(\ell)} \mathscr{U}_2^{(\ell)}...\mathscr{U}_{n-1}^{(\ell)} \mathscr{U}_n^{(\ell)}$
\end{algorithmic}
\label{main_algo}
\end{algorithm}

The convergence criterion for the proposed TCAM-TT algorithm is defined {via a threshold on the relative change, say $\epsilon$, in the successive estimation of the factors},
\begin{equation}\label{eq: criterion}
\epsilon =\sum_i \frac{ \|\mathscr{U}_i^{(\ell+1)} -\mathscr{U}_i^{(\ell)} \|_F}{ \| {\mathscr{U}}^{(\ell)}_i\|_F}
\end{equation}
where $\ell$ is the iteration parameter and $maxiter$ is maximum iterations. {The algorithm will stops either when $\ell$ reaches $maxiter$ or when the relative error $\epsilon \leq tot$ for some predefined tolerance parameter $tot$.}

\section{Complexity Analysis}\label{compl}
%The computation complexity mainly includes three parts, Least Square Calculation Complexity (LSCC), Tensor Connect Product Reconstruction Complexity (TCPRC).

The main algorithm in the computation is the computation of the least squares. This computation is performed for each slice of each tensor train factors $\mathscr{U}_{k=1,\cdots, n}(:, i_k, :), i_k\in 1,\cdots I_k$. The matrix corresponding to the  least square problem $(\|{\bf A}x -{\bf b}\|_F)$ satisfies  ${\bf A} \in \mathbb{R}^{P_{k,i_k} \times (r_{k-1}r_k)}$, where $P_{k, i_k}$ is the number of observed entries in the $i_k^{\text{th}}$ in $\mathscr{X}_{[k]}$. For the analysis, we assume that  all ranks are the same, or $r_k=r$. Since the complexity of pseudo-inverse of $d_1\times d_2$ matrix is $O(d_1^2d_2)$ \cite{Complexity}, this complexity is given as $O(P_{k,i_k} r^2_{k-1}r^2_k)$. Thus, the overall complexity in each iteration is given as $O(nPr^4)$, where  $P$ is the total number of observed entries.

%Least square calculation complexity is the complexity for solving each slice of tensor train factors, e.g. $\mathscr{U}_{k=1,\cdots, n}(:, i_k, :), i_k\in 1,\cdots I_k$, where the matrix for this least square problem $(\|{\bf A}x -{\bf b}\|_F)$,   ${\bf A} \in \mathbb{R}^{P_{k,i_k} \times (r_{k-1}r_k)}$, where $P_{k, i_k}$ is the number of observed entries in the $i_k^{\text{th}}$ in $\mathscr{X}_{[k]}$. The pseudo-inverse of ${\bf A}$ takes $O(P_{k,i_k} r^2_{k-1}r^2_k) +O(r^3_{k-1}r^3_k)$, where $O(P_{k,i_k} r^2_{k-1}r^2_k)$ dominates.
%If assuming that all TT-Rank are the same, denoted as $r$,then the least square calculation complexity for $\mathscr{U}_k$ is $O(Pr^4)$, where $P$ is the total number of observed entries, indicating that the LSCC for updating all $\mathscr{U}_k$ has the complexity $O(nPr^4)$.

%Tensor Connect Product Reconstruction Complexity is the complexity for calculating the connected product of all the tensor train decomposition terns besides the one to be updated. As described in Fig \ref{ConProd}, where the connect product between $\mathscr{U}_j \mathscr{U}_{j+1}$ is $O(r_{j-1}I_j \times r_j \times r_{j+1} \times I_{j+1})$. 
%If assuming that all TT-Rank are the same, denoted as $r$,  and each mode has $I$ entries,  , the complexity for tensor connect product is $O(I^{n-1}r^3)$

\section{Numerical Results}
\label{sec:5}
In this section, we compare our proposed TCAM-TT algorithm with Tensor Completion by alternating Minimization after Tensor Mode Matrization  (TCAM-TMM),  SiLRTC-TT algorithm as proposed in \cite{phien2016efficient} and tSVD algorithm as proposed in \cite{zhang2013novel}. We briefly describe these algorithms below.

\subsection{TCAM-TMM}
The first tensor completion algorithm is {TCAM-TMM} algorithm where a tensor is unfolded into matrix and alternating minimization \cite{jain2013low} is used to solve the resulting matrix completion problem. For an order-$n$ tensor with a given set of tensor train rank, TCAM-TMM algorithm uses the $n-1$ possible ways of of tensor mode matricization (
%\textcolor{blue}{the reader is again referred back to the definition used in this paper - }
see section \ref{sec:2}, Definition \ref{def2}. 

In particular, let $\mathscr{X}_{\floor{k:k=1,\cdots, n-1}} \in \mathbb{R}^{(\prod_{t=1}^k I_t) \times (\prod_{t=k+1}^N I_t)}$ be the tensor the $k^{\text{th}}$ mode matrization of tensor $\mathscr{X} \in \mathbb{R}^{I_1 \times ... \times I_n }$, thus TCAM-TMM solves the tensor completion problem by solving the following matrix completion problem,
\begin{equation}
\min_{{\bf U}_k, {\bf V}_k} \|  \mathscr{P}_{\Omega_{\floor{k}}} \circ({\bf U}_k {\bf V}_k^\top -\mathscr{X}_{\floor{k}})\|_F,
\end{equation}
where $ \mathscr{P}_{\Omega_{\floor{k}}} $ is the binary tensor after the $k^{\text{th}}$ tensor mode matrization, ${\bf U}_k$ and ${\bf V}_k$ are the low-rank
%(\textcolor{red}{$\leftarrow$, Wenqi: Here should we specify the rank, that I believe will come from the TT rank}) 
factorization terms of the tensor $\mathscr{X}$ after the $k^{\text{th}}$ tensor mode matrization and the rank is the $r_k$, which is selected from the tensor train rank.

\subsection{SiLRTC-TT}
The second tensor completion algorithm is SiLRTC-TT algorithm as proposed in \cite{phien2016efficient}, which completes tensors by taking orders to do matrix completion after tensor mode unfolding and recovery the tensor by weighted summarization of the tensor after each matrix completion. It is selected as it has been shown to have the best performance in \cite{phien2016efficient}.

\subsection{tSVD}
The third tensor completion algorithm is the tubal-SVD (t-SVD) based algorithm as proposed in \cite{zhang2013novel,Zhang_TSP}. This algorithm works by minimizing the nuclear norm of a block circulant matrix that is formed out of the slices of the tensor. This algorithm is selected as it shows very good performance for video completion.\\

The performance of all these algorithms are measured by the Recovery Error at Missing Entries (REME), defined as
\begin{equation}\label{eq: measure}
\text{REME} = \frac{\| \mathscr{X}_{-\Omega} - \mathscr{X}^{(R)}_{-\Omega}\|_F}{\| \mathscr{X}_{-\Omega} \|_F},
\end{equation}
where $\mathscr{X}_{-\Omega}$ represents missing entries in the original tensor, $ \mathscr{X}^{(R)}_{-\Omega}$ represents missing entries in the recovered tensor.

%*********************************************  Synthetic Data ************************************************************************
\subsection{Synthetic Data}
In this section, we consider a completion problem of a $4$ dimensional tensor $\mathscr{X} \in \mathbb{R}^{20 \times 20 \times 20\times 20}$ with TT-Rank $[1, 5, 100, 5, 1]$ without loss of generality. The tensor is generated by a sequence of connected tensors $\mathscr{U}_{i:i=1,2,3,4}$, and all the  entries in $\mathscr{U}_i$ are sampled from independent standard normal distribution.
 
The 4-D tensor $\mathscr{X}$ with a pre-defined tensor train rank has $3$ tensor mode matrization, and each tensor mode matrixzation, denoted as $\mathscr{X}_{\floor{k}, k=1,2,3}$ generates a matrix completion problem of $\mathscr{X}_{\floor{1}} \in \mathbb{R}^{20 \times 8000}$, $\mathscr{X}_{\floor{2}} \in \mathbb{R}^{400 \times 400}$ and $\mathscr{X}_{\floor{3}} \in \mathbb{R}^{8000 \times 20}$ with rank $5$, $100$ and $5$ respectively.    The completion results after each tensor mode matrization (TMM) are denoted as TCAM-TMM1, TCAM-TMM2, and TCAM-TMM3.

The error tolerance $tot$ for all algorithm is set to be $10^{-4}$. Thus any REME that is lower than $10^{-4}$ is regarded as a perfect completion.
The maximum iteration, $maxiter$, is set to be $100$ for TCAM-TT  and $1000$ for TCAM-TMM1, TCAM-TMM2, TCAM-TMM3, SiLRTC-TT, and tSVD algorithm. 

\begin{figure}[ht]
\includegraphics [trim=0.5in 2.6in 1in 2.9in, keepaspectratio, width=0.4\textwidth] {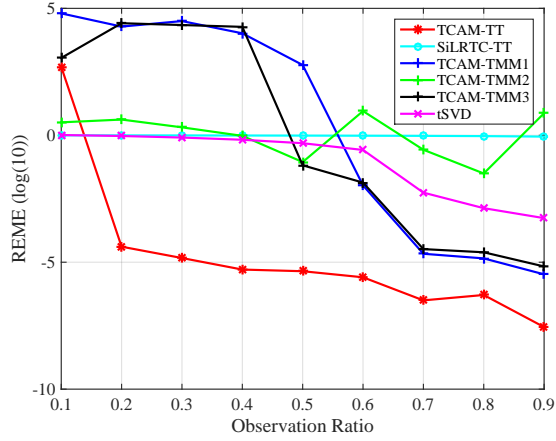}
\centering
\vspace{-.1in}
\caption{Tensor completion under tensor train model for synthetic tensor $\mathscr{X} \in \mathbb{R}^{20 \times 20 \times 20 \times 20}$ with TT-rank $[1, 5, 100, 5,1]$}
\label{synthethic4}
\end{figure}

The simulation in Fig \ref{synthethic4} shows the REME at log10 scale for observation ratio from $10\%$ to $90\%$ of all algorithms and each plotted point is the average of 12 independent repeated experiment. 
TCAM-TT algorithm performs the best as it achieves perfectly tensor recovery for any observation ratio higher $20\%$ while 
TCAM-TMM1, TCAM-TMM3, and tSVD achieve perfectly recovery at the sampling ratio $70\%$, $70\%$, and $90\%$ . 
TCAM-TMM2 and SiLRTC-TT are not effective in tensor completion in this case as the recovery errors for the two algorithms are around $1$.
TCAM-TT algorithm achieve the best performances as it consider all the tensor train rank together and the updating of each alternating minimization step maintains the TT-Rank property. In contrast, TACM-TMM, SiLRTC-TT and tSVD algorithm consider each tensor train rank independently, thus the completion results fit one specific rank of the TT-Rank very well but may not fit all the TT-Rank, which leads to the lower performance.

\begin{figure}[ht]
\includegraphics [trim=0.5in 2.6in 1in 2.9in, keepaspectratio, width=0.4\textwidth] {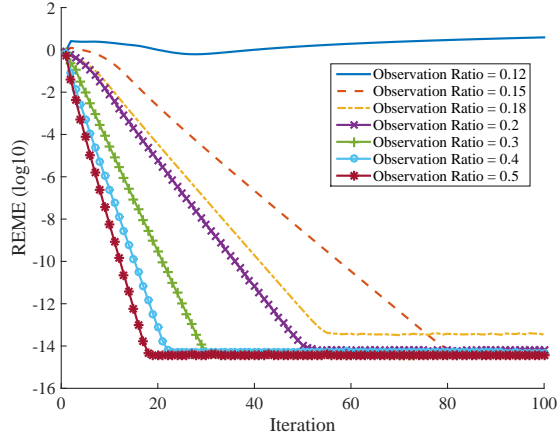}
\centering
\vspace{-.1in}
\caption{TCAM-TT algorithm convergence for synthetic tensor $\mathscr{X} \in \mathbb{R}^{20 \times 20 \times 20 \times 20}$ with TT-rank $[1, 5, 100, 5, 1]$}
\label{syn4_convergence}
\end{figure}

In addition to better recovery as compared with the other algorithms, TCAM-TT algorithm converges within less number of alternating minimization iterations. The convergence performance of TCAM-TT algorithm for observation ratio from 12\% to 50\% is shown in Fig \ref{syn4_convergence}. We note that for any observation ratio larger than $15\%$ TCAM-TT takes less than $80$ iterations. Typically, the larger the observation ratio, the less iteration it takes for TCAM-TT to converge. For example, the tensor with $50\%$ observation ratio only takes $19$ iterations to converge with the REME being $10^{-16}$. {This fast convergence is in part  due to good initialization when more data is available.}

%%%%%%%%%%%%%%%%%%%%%%%%%%  YaleFace  %%%%%%%%%%%%%%%%%%%%%%%%%%%%%%%%%%%%%%%%
\subsection{Extended YaleFace Dataset B}
Extended YaleFace Dataset B \cite{GeBeKr01} is a dataset that includes 38 people with 9 poses under 64 illumination conditions. Each  image has the size of $192 \times 168$, where we down-sample the size of each image to $48 \times 42$ for ease of calculation. We consider the images for 38 people under 64 illumination within 1 pose by reshaping the data into a tensor $\mathscr{X} \in \mathbb{R}^{48 \times 42\times 64 \times 38}$.
TT-rank is estimated to be $ [1, 31,  137,  31, 1]$, which gives $10\%$ error for fitting the dataset when there are no missing entries. 
Missing entries are sampled by assuming that data is entry-wise missing with probability $p$, where $p$ changes from $10\%$ to $90\%$.
The error tolerance $tot$ for all algorithm is again set to be $10^{-4}$ and 
the maximum iteration $maxiter$ is set to be $20$ for TCAM-TT while the maximum iteration for TCAM-TMM1, TCAM-TMM2, TCAM-TMM3, SiLRTC-TT and tSVD are all set to be $1000$. 

\begin{figure}[h]
\includegraphics [trim=0.5in 0.6in 0.5in 0.5in, keepaspectratio, width=0.45\textwidth] {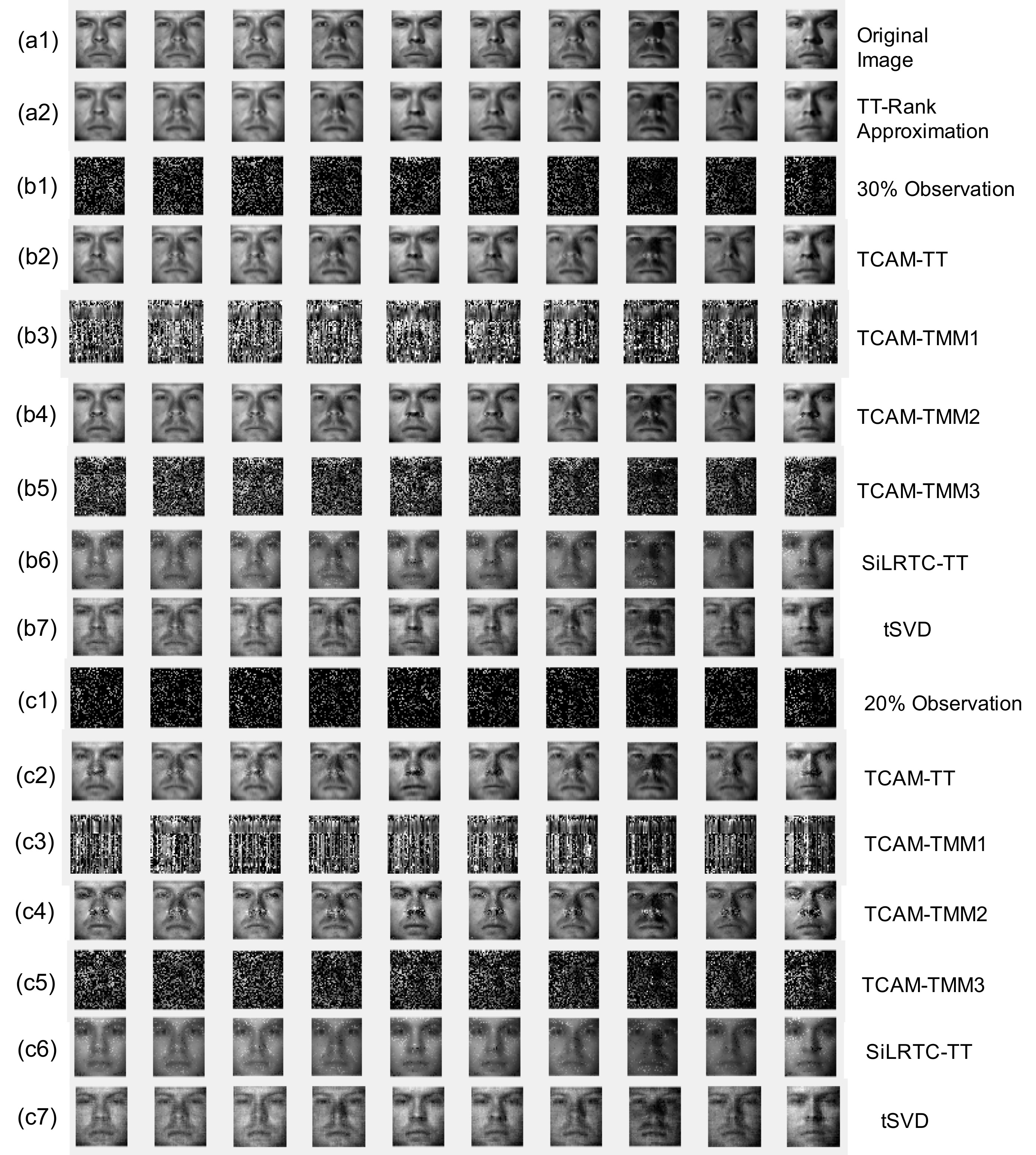}
\centering
\vspace{0.1in}
\caption{Tensor completion for Extended YaleFace Dataset B. 
(a1) Original 10 images selected from the dataset. 
(a2) 10 images after TT-rank $[1, 31, 137, 31, 1]$ approximation when no missing entries where the recovery error is $10\%$.
(b1-b7) shows missing data, TCAM-TT completed data, TCAM-TMM1 completed data, TCAM-TMM2 completed data,TCAM-TMM3 completed data, SiLRTC-TT completed data and tSVD completed under the scenario when each entry is missing with probability $0.7$.  
%The completion error is $14.64\%$(TCAM-TT), $281.4\%$(TCAM-TMM1), $21.7\%$(TCAM-TMM2), $147.6\%$(TCAM-TMM3), $31.2\%$(SiLRTC-TT) and $18.6\%$(tSVD).
(c1-c7) shows missing data, TCAM-TT completed data, TCAM-TMM1 completed data, TCAM-TMM2 completed data,TCAM-TMM3 completed data, SiLRTC-TT completed data and tSVD completed under the scenario when each entry is missing with probability $0.8$.  
%The completion error is $22.31\%$(TCAM-TT), $204.2\%$(TCAM-TMM1), $72.91\%$(TCAM-TMM2), $117.22\%$(TCAM-TMM3), $33.97\%$(SiLRTC-TT) and $21.85\%$(tSVD).
}
\label{YaleFace}
\end{figure}

\begin{figure}[h]
\includegraphics [trim=0.5in 0.6in 0.5in 0.5in, keepaspectratio, width=0.4\textwidth] {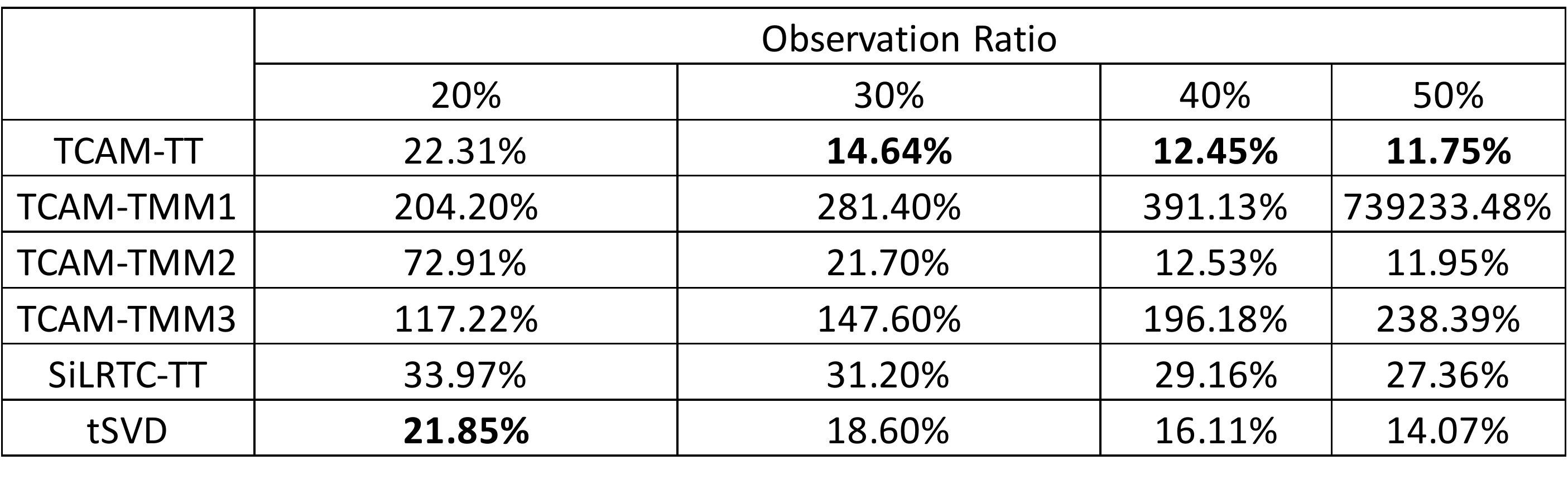}
\centering
\vspace{0.1in}
\caption{ REME for Extended YaleFace Dataset B Completion at Observation Ratio $20\%$, $30\%$, $40\%$ and $50\%$}
\label{YaleFaceTable}
\end{figure}

The simulation results shown in Fig \ref{YaleFace} describe the completed images under $70\%$ and $80\%$ observation and the table in Fig \ref{YaleFaceTable} shows the REME values for each algorithm.

\begin{figure}[h]
\includegraphics [trim=0.5in 2.6in 1in 2.9in, keepaspectratio, width=0.4\textwidth] {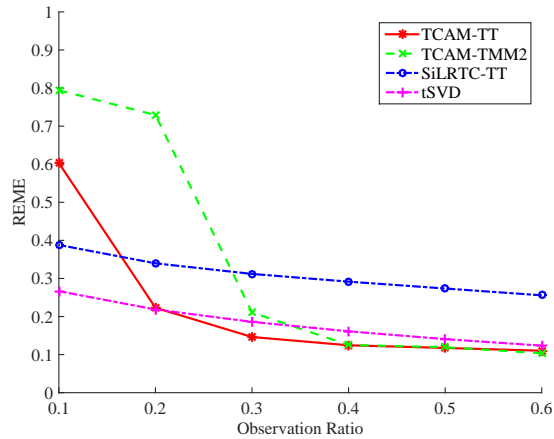}
\centering
\caption{ REME versus observation ratio from 10\% to 60\% for Extended YaleFace Dataset B}
\label{YaleFace_Plot}
\end{figure}

Noting that under $70\%$ observation ratio, only TCAM-TT, TCAM-TMM2, SiLRTC-TT and tSVD can complete the images while TCAM-TMM1 and TCAM-TMM3 algorithm can not, 
thus only TCAM-TT, TCAM-TMM2, SiLRTC-TT and tSVD algorithm are considered for investigating the relation between observation ratio and REME for Extended YaleFace Dataset B, as shown in Fig \ref{YaleFace_Plot}. 

SiLRTC-TT and tSVD algorithm both show stable recovery result for all sampling ratios. When sampling ratio decreases from $60\%$ to $10\%$, the recovery error increases from $25.57\%$ and $12.3\%$ to $38.80\%$ and $26.7\%$, which shows the stable performance in all the algorithm. However, the stable performance comes with the cost of blurry recovery, as shown in Fig \ref{YaleFace} (b6) and (b7), where although the recovery result is smooth, each image is less sharp in resolution. tSVD algorithm performs better than SiLRTC-TT algorithm under any observation ratio.

Both TCAM-TT and TCAM-TMM2 algorithm have shown good recovery when the sampling ratio is greater than $40\%$ and the increasing of error of recovery when the sampling ratio becomes lower. 
The recovery result for TCAM-TMM2 starts to degrade at $40\%$ sampling ratio and the error increases faster thanTCAM-TT algorithm, as TCAM-TMM2 does not capture the tensor structure in first and third unfolding of the tensor.
TCAM-TT algorithm shows the best result for all sampling ratio larger than $20\%$. 

%%%%%%%%%%%%%%%%%%%%%%%%%%%%%%%%%%  Video Data %%%%%%%%%%%%%%%%%%%%%%%%%%%%%%%%%%%
\subsection{Video Data}
The Video data we used is high speed camera video for bullet we downloaded from Youtube \cite{Video_Web} with $85$ frames in total and each frame is consisted by a $100 \times 260 \times 3$ color image.
The video data is regarded as a 4-mode tensor $\mathscr{X} \in \mathbb{R}^{100 \times 260 \times 3 \times 85}$. 
Different from the tensor constructed from Extended YaleFace Dataset B where the $4^{\text{th}}$ mode of the tensor that represents different persons only has weak connections, the $4^{\text{th}}$ mode of the tensor built from the video data owns a stronger connection and the lower rank property is more likely to hold as the $4^{\text{th}}$ mode of the tensor represents the time series and any frame is easily to be represented by the linear combination of its previous frame and its next frame. 
The TT-rank is estimated to be $[1, 29, 99, 19, 1]$, which gives $6.33\%$ error for fitting the dataset when there are no missing data. 

\begin{figure}[h]
\includegraphics [trim=0.5in 0.6in 0.5in 0.5in, keepaspectratio, width=0.45\textwidth] {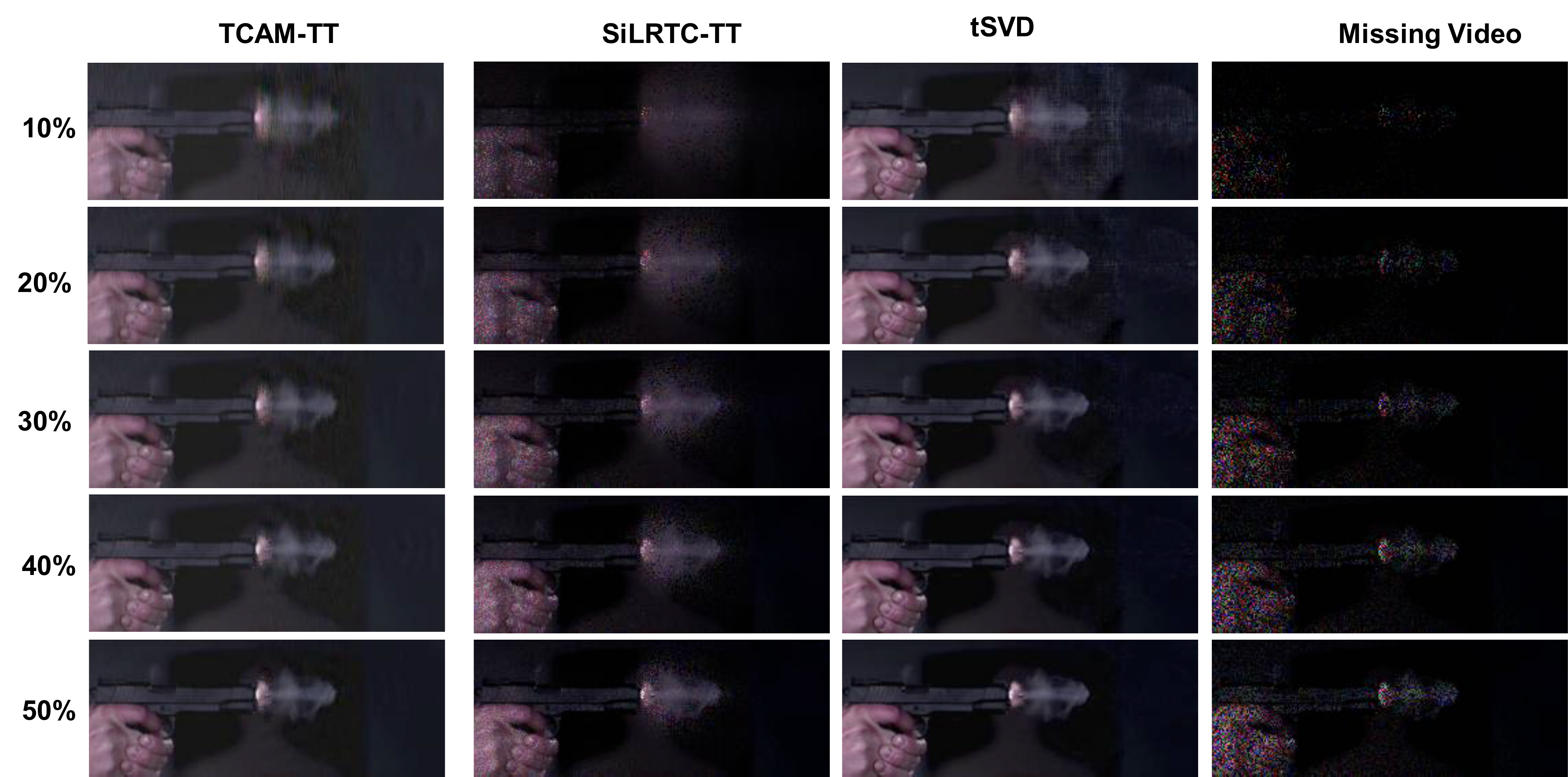}
\centering
\vspace{0.1in}
\caption{Video completion by TCAM-TT, SiLRTC, and tSVD under observation ratio $10\%$, $20\%$, $30\%$, $40\%$, and $50\%$. The $1^{st}$ frame of the video is displayed. }
\label{Video}
\end{figure}

This video is selected as under high speed, gun and hand are almost still while smoke and bullet are movable, which could show the algorithm recovery performance on both still and dynamic objects within video. The $1^{\text{st}}$ frame of the recovered video image is shown in Fig \ref{Video}, where entries in video are set to be missing independently with probability $p$,  which changes from $10\%$ to $90\%$ at the step of $10\%$.

TCAM-TT, TCAM-TMM1, TCAM-TMM2, TCAM-TMM3, SiLRTC-TT and tSVD algorithm are implemented while only TCAM-TT, SiLRTC-TT and tSVD algorithm are able to complete the video when the observation ratio is less than $70\%$, showing the advantage of tensor completion as compared with matrix completion when high order data is considered. 
Completion result in Fig \ref{Video} shows that TCAM-TT algorithm out performs than SiLRTC-TT algorithm as the still objects can be recovered completely and dynamic objects can be recovered smoothly for adjacent frames. In contrast, SiLRTC-TT algorithm does not complete the video completely as the hand in the figure is blurred by the uncompleted dots. 

\begin{figure}[h]
\includegraphics [trim=0.5in 2.6in 1in 2.9in, keepaspectratio, width=0.4\textwidth] {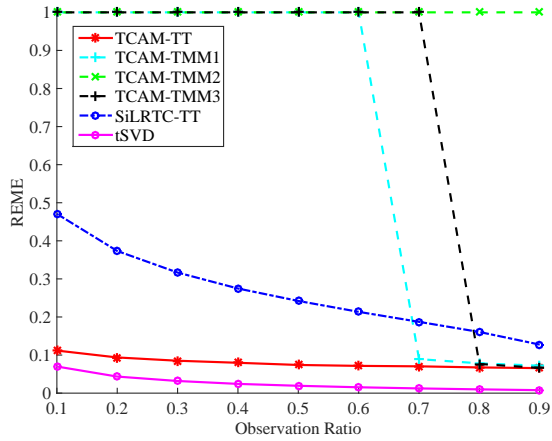}
\centering
\caption{ REME versus observation ratio from 10\% to 90\% for Video Completion}
\label{Video_Plot}
\end{figure}

Figure \ref{Video_Plot} shows the REME for observation from $10\%$ to $90\%$, where any REME larger than $1$ is set to be $1$ for ease of visualization. 
The completion results for SiLRTC-TT degrades faster than TCAM-TT algorithm when the observation ratio decreases since the video becomes more uniformly dark and more blurry when the missing ratio increases from $50\%$ to $90\%$. 
In the video completion, tSVD performs the best of all the algorithm, which benefits from the advantage of Fourier transform that is applied on time series. The error in the proposed algorithm is limited by the error in TT-rank approximation of the actual data with the chosen rank. The proposed algorithm performs the best among the other algorithms, and in specific as compared to the other algorithms that exploit the TT-rank structures and the matrix unfolding based approaches.

%%%%%%%%%%%%%%%%%%%%%%%%%%%%%%%%%%%%%   Seismic Data   %%%%%%%%%%%%%%%%%%%%%%%%%%%%%%%%%%%%%%
\subsection{Seismic Data}

In thie subsection, we widh to complete pre-stack seismic
records from incomplete spatial measurements. The pre-stack seismic data can be viewed as a 5D data or a
fifth order tensor consisting of one time or frequency dimension
and four spatial dimensions describing the location of the
detector and the receiver in a two dimensional plane. This data
can then be described in terms of the original $(r_x,r_y,s_x,s_y)$ coordinate
frames or in terms of midpoint receivers and offsets
$(x, y,h_x,h_y)$ \cite{Seismic}.  We use the dataset from  \cite{Seismic}, where the sources and receivers are placed on  a $16\times 16$ grid with $50m$ shot forming a tensor $\mathscr{X} \in \mathbb{R}^{16 \times 16 \times 17 \times 17 \times 150}$. We approximate the TT- Rank of this tensor as $[1, 1, 1, 9, 41, 1]$.

%Seismic data is a 5-D tensor in which sources and receivers are placed on  a $16\times 16$ grid with $50m$ shot and receiver spacing in the $x$ and $y$ directions. 

%The synthetic seismic data is a tensor $\mathscr{X} \in \mathbb{R}^{16 \times 16 \times 17 \times 17 \times 150}$ where the approximated TT- Rank is $[1, 1, 1, 9, 41, 1]$. 

%In the scenario of seismic data, data is time series - wise missing rather than entry-wise missing as happened in Extended YaleFace Dataset B and video data. 

\begin{figure}[h]
\includegraphics [trim=0.5in 2.6in 1in 2.9in, keepaspectratio, width=0.4\textwidth] {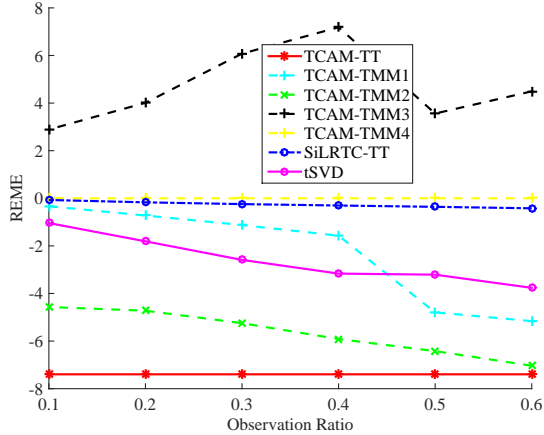}
\centering
\caption{ REME versus observation ratio from 10\% to 60\%  for Seismic Data}
\label{Seimic_Plot}
\end{figure}

The results in Figure \ref{Seimic_Plot} illustartes that the proposed algorithm performs significantly better than other compared algorithms for different observation ratios from 10\% to 60\%.

\section{Conclusion}
\label{sec:6}
We proposed a novel algorithm for data completion using tensor train decomposition. Unlike the current methods exploiting this format, our algorithm exploits the matrix product state representation and uses alternating minimization over the low rank factors for completion. As a future work we will derive provable performance guarantees on tensor completion using the proposed algorithm. In this context, the statistical machinery for proving analogous results for the matrix case \cite{xx} can be used. We will also look at parallelizing this algorithm and make it more efficient in terms of implementation, especially when forming and storing the intermediate tensors from the estimated factors. 

\section{appendix}
%================================= Lemma 1 ======================================
\subsection{Proof of Lemma 1}\label{proof0}
\proof
Let ${\bf M} ={\bf M}_1 {\bf M}_2$, thus 
\begin{equation}
{\bf M}(j_1, j_2) =\sum_{j=1}^{r_1} {\bf M}_1(j_1,j) {\bf M}_2(j, j_2)
\end{equation}
where ${\bf M}(j_1, j_2)$ locates at $\text{vec}({\bf M}_1{\bf M}_2)(j_1 + (j_2-1)I_1, 1)$.

Let ${\bf T}_1 \in \mathbb{R}^{(I_1I_2) \times (r_1I_2)} = {\bf I}^{(I_{2})} \otimes {\bf L}(\mathscr{M}_1)$ and ${\bf T}_2\in \mathbb{R}^{(r_1I_{2})  \times 1} = {\bf L}(\mathscr{M}_2)$, and ${\bf T} \in\mathbb{R}^{I_1I_2 \times 1} = {\bf T}_1{\bf T}_2$, thus 
\begin{equation}
\begin{split}
&{\bf T}(j_1+(j_2-1)I_1, 1) \\
= & \sum_{j=1}^{r_1I_2}{\bf T}_1(j_1+(j_2-1)I_1, j) {\bf T}_2(j, 1)\\
= & \sum_{j=(j_2-1)r_1+1}^{j_2r_1} {\bf T}_1(j_1+(j_2-1)I_1, j) {\bf T}_2(j, 1)\\
= & \sum_{j=1}^{r_1}{\bf M}(j_1, j){\bf M}_2(j, j_2)
\end{split}
\end{equation}
We conclude that any ${j_1 + (j_2-1)I_1}^\text{th}$ entry on the left hand side is the same as that on the right hand side, thus we prove our claim. 
\endproof
%================================ Lemma 2 ===============================
\subsection{Proof of Lemma 2}\label{proof1}
\proof
Based on definition of tensor permutation in \eqref{eq: TensorPermutation}, on the left hand side, the $(j_1,...., j_n)$ entry of the tensor is 
\begin{equation} \label{eq: prf1}
\mathscr{X}^{P_i}(j_1,...,j_n) = \mathscr{X}(j_{n-i+2},..., j_{n}, j_{1},..., j_{n-i+1}).
\end{equation}

On the right hand side, the $(j_1,...., j_n)$ entry of the tensor gives
\begin{equation}
\begin{split}
& f(\mathscr{U}_i\cdots \mathscr{U}_{i-1})(j_1,\cdots, j_n)\\
=&\text{Trace}(\mathscr{U}_i(:, j_1,:)\mathscr{U}_{i+1}(:, j_2,:)... \mathscr{U}_n(:, j_{n-i+1},:)\\
& \mathscr{U}_1(:, j_{n-i+2}, :) \cdots \mathscr{U}_{i-1}(:, j_n, 1)).
\end{split}
\end{equation}
Since trace is invariant under cyclic permutations, we have
\begin{equation}
\begin{split}
&\text{Trace}(\mathscr{U}_i(:, j_1,:)\mathscr{U}_{i+1}(:, j_2,:)... \mathscr{U}_n(:, j_{n-i+1},:)\\
& \mathscr{U}_1(:, j_{n-i+2}, :) \cdots \mathscr{U}_{i-1}(:, j_n, 1))\\
=&\text{Trace}(
 \mathscr{U}_1(:, j_{n-i+2}, :) \cdots \mathscr{U}_{i-1}(:, j_n, 1)\\
&\mathscr{U}_i(:, j_1,:)\mathscr{U}_{i+1}(:, j_2,:)... \mathscr{U}_n(:, j_{n-i+1},:))\\
= &f(\mathscr{U}_1\cdots \mathscr{U}_{n})(j_{n-i+2},\cdots, j_n, j_1,\cdots, j_{n-i+1}),
\end{split}
\end{equation}
which equals to the right hand side of equation \eqref{eq: prf1} based on \eqref{eq: eq6}. 
Since any entries in $\mathscr{X}^{P_i}$ are the same as those in $\mathscr{U}_i \mathscr{U}_{i+1} \cdots \mathscr{U}_n \mathscr{U}_1 \cdots \mathscr{U}_{i-1}$, the claim is proved.
\endproof

%===============================Lemma 3=========================================
\subsection{Proof of Lemma 3} \label{proof2}
\proof
First we note that tensor permutation does not change tensor Frobenius norm as all the entries remain the same as those before the permutation. 
Thus, when $i\neq 1$, we permute the tensor inside the Frobenius norm in \eqref{eq: tran1} and get the equivalent equation  as
\begin{equation}\label{eq: p0}
\mathscr{U}_i = \argmin_{\mathscr{Y}} \| \mathscr{P}^{P_i}_\Omega \circ (f(\mathscr{U}_1 \cdots\mathscr{U}_{i-1}\mathscr{Y} \mathscr{U}_{i+1}\cdots \mathscr{U}_n))^{P_i} -\mathscr{X}^{P_i}_\Omega \|_F^2.
\end{equation}

Based on Lemma \ref{lemma1}, we have
\begin{equation}
(f(\mathscr{U}_1 \cdots\mathscr{U}_{i-1}\mathscr{Y} \mathscr{U}_{i+1}\cdots \mathscr{U}_n))^{P_i} =f(\mathscr{Y} \mathscr{U}_{i+1}\cdots \mathscr{U}_n \mathscr{U}_1 \cdots\mathscr{U}_{i-1}),
\end{equation}
thus equation \eqref{eq: p0} becomes 
\begin{equation}\label{eq: p0_1}
\mathscr{U}_i = \argmin_{\mathscr{Y}} \| \mathscr{P}^{P_i}_\Omega \circ f(\mathscr{Y} \mathscr{U}_{i+1}\cdots \mathscr{U}_n \mathscr{U}_1 \cdots\mathscr{U}_{i-1}) -\mathscr{X}^{P_i}_\Omega \|_F^2.
\end{equation}
Comparing \eqref{eq: p0_1} and \eqref{eq: T3}, we have $\mathscr{P}_\Omega,\mathscr{X}_\Omega$ and $\mathscr{U}_2\cdots\mathscr{U}_n $ in \eqref{eq: T3} become $\mathscr{P}_\Omega^{\top_i}, \mathscr{X}_\Omega^{\top_i}$ and $\mathscr{U}_{i+1} \cdots \mathscr{U}_n\mathscr{U}_1 \cdots \mathscr{U}_{i-1}$ in\eqref{eq: p0_1} respectively. Thus we prove our claim. 
\endproof

%===================================================================================
\subsection{Proof of Lemma 4} \label{proof3}
\proof
\begin{equation}
\begin{split}
\text{Trace}(A \times B) &= \sum_i^{r_1}\left(\sum_j^{r_2} {\bf A}(i,j) {\bf B}(j, i ) \right) \\
&= \sum_i^{r_1}\sum_j^{r_2}  {\bf A}(i,j) {\bf B}^\top(i, j ) \\
& = vec({\bf A})^\top vec({\bf B}^\top)
\end{split}
\end{equation}
\endproof

%==== Bib files and style =======
%\vspace{-3mm}
\bibliographystyle{IEEEbib}
\bibliography{Ref,CAMSAP_2015_bib,bibtensor}

\end{document}